%&amstex
\input amstex
\input amsppt.sty
\magnification=\magstep1
\hsize=30truecc
\vsize=22.2truecm
\baselineskip=16truept
\NoBlackBoxes
\TagsOnRight \pageno=1 \nologo
\def\Z{\Bbb Z}
\def\N{\Bbb N}

\def\l{\left}
\def\r{\right}
\def\bg{\bigg}
\def\({\bg(}
\def\[{\bg\lfloor}
\def\){\bg)}
\def\]{\bg\rfloor}
\def\t{\text}
\def\f{\frac}

\def\bi{\binom}
\def\eq{\equiv}

\def\ls{\leqslant}
\def\gs{\geqslant}
\def\mo{\roman{mod}}

\def\jacob #1#2{\left(\frac{#1}{#2}\right)}
\def\qbinom #1#2#3{\left[\matrix {#1} \\ {#2}\endmatrix\right]_{#3}}
\def\pmod #1{\ (\roman{mod}\ #1)}

\def\Mm#1#2#3{\thickfracwithdelims[]\thickness0{#1}{#2}_{#3}}

\def\Proof{\noindent{\it Proof}}

\def\Remark{\medskip\noindent{\it  Remark}}

\def\Ack{\medskip\noindent {\bf Acknowledgment}}
\hbox {Preprint, {\tt arXiv:0910.4170}}
\bigskip
\topmatter
\title Some $q$-congruences related to $3$-adic valuations\endtitle
\author Hao Pan and Zhi-Wei Sun\endauthor
\leftheadtext{Hao Pan and Zhi-Wei Sun}
\affil Department of Mathematics, Nanjing University\\
 Nanjing 210093, People's Republic of China
\\{\tt haopan79\@yahoo.com.cn,
zwsun\@nju.edu.cn}
\endaffil
\abstract In 1992 Strauss, Shallit and  Zagier proved that for any positive integer $a$ we have
$$\sum_{k=0}^{3^a-1}\binom{2k}{k}\eq0\ (\mo\ 3^{2a})$$ and furthermore
$$\f1{3^{2a}}\sum_{k=0}^{3^a-1}\binom{2k}k\eq1 (\mo\ 3).$$ Recently a $q$-analogue of
the former congruence was conjectured by Guo and Zeng. In this paper we prove
the conjecture of Guo and Zeng, and also give a $q$-analogue of the latter congruence.
\endabstract
\thanks 2010 {\it Mathematics Subject Classification}.\,Primary 11B65;
Secondary 05A10, 05A30, 11A07, 11S99.
\newline\indent The second author is the corresponding author.
He was supported by the National Natural Science
Foundation (grant 10871087) and the Overseas Cooperation Fund (grant 10928101) of China.
\endthanks
\endtopmatter
\document

\heading{1. Introduction}\endheading

Partially motivated by the work of Pan and Sun [PS], Sun and Tauraso
[ST1] proved that for any prime $p$ and $a\in\Z^+=\{1,2,3,\ldots\}$
we have
$$\sum_{k=0}^{p^a-1}\bi{2k}k\eq\l(\f{p^a}3\r)\ (\mo\ p^2),$$
where $(-)$ denotes the Legendre symbol. (See also [ST2, ZPS,
S09a, S09b, S09c] for related results.) When checking whether there are
composite numbers $n$ such that
$$\sum_{k=0}^{n-1}\bi{2k}k\eq\l(\f n3\r)\ (\mo\ n^2),$$
Sun and Tauraso found that
$$\nu_3\(\sum_{k=0}^{3^a-1}\bi{2k}k\)\gs 2a\quad\ \t{for}\ a=1,2,3,\ldots,\tag1.1$$
where $\nu_3(m)$ denotes the $3$-adic valuation of an integer $m$
(i.e., $\nu_3(m)=\sup\{a\in\N:\ 3^a\mid m\}$ with
$\N=\{0,1,2,\ldots\}$). However, a refinement of this was proved earlier by
Strauss, Shallit and  Zagier [SSZ] in 1992.

\proclaim{Theorem 1.1 {\rm (Strauss, Shallit and Zagier [SSZ])}} For any
$a\in\Z^+$ we have
$$\sum_{k=0}^{3^a-1}\bi{2k}k\eq 3^{2a}\ (\mo\ 3^{2a+1}).\tag1.2$$
Furthermore,
$$\f{\sum_{k=0}^{n-1}\bi{2k}k}{n^2\bi{2n}n}\eq-1\ (\mo\ 3)\quad\t{for all}\ n\in\Z^+.$$
\endproclaim

Recall that the usual $q$-analogue of $n\in\N$ is
$$[n]_q=\f{1-q^n}{1-q}=\sum_{0\ls k<n}q^k$$
which tends to $n$ as $q\to1$. For $d\in\Z^+$ the $d$-th cyclotomic
polynomial is given by
$$\Phi_d(q)=\prod\Sb r=1\\(r,d)=1\endSb^d\l(q-e^{2\pi ir/d}\r)\in\Z[q].$$
Given a positive integer $n>1$ we obviously have
$$[n]_q=\f{q^n-1}{q-1}=\prod_{k=1}^{n-1}\l(q-e^{2\pi i k/n}\r)=\prod\Sb d\mid n\\d>1\endSb \Phi_d(x).$$
It is well known that if $d_1,d_2\in\Z^+$ are distinct then $\Phi_{d_1}(q)$ and $\Phi_{d_2}(q)$
are relatively prime in the polynomial ring
$\Z[q]$. If $p$ is a prime and $a$ is a positive integer, then
$$\Phi_{p^a}(q)=\f{q^{p^a}-1}{q^{p^{a-1}}-1}=[p]_{q^{p^{a-1}}}\ \ \t{and}\ \ [p^a]_q=\prod_{j=1}^a\Phi_{p^j}(q).$$

For $n,k\in\N$ the usual $q$-analogue of the binomial coefficient
$\bi nk$ is the folloqing $q$-binomial coefficient:
$$\Mm nkq=\cases([n]_q\cdots[n-k+1]_q)/([1]_q\cdots[k]_q)&\t{if}\ 0<k\ls n,
\\1&\t{if}\ k=0,\\0&\t{if}\ k>n.\endcases$$
Note that $\Mm nkq\to \bi nk$ as $q\to1$. Many combinatorial
identities and congruences involving binomial coefficients have
their $q$-analogues (cf. [St]).

Recently Guo and Zeng [GZ] proposed a $q$-analogue of (1.1), namely
they formulated the following conjecture.

\proclaim{Conjecture 1.2 {\rm (Guo and Zeng [GZ,
Conjecture 3.5])}} Let $a$ be a positive integer. Then
$$\sum_{k=0}^{3^am-1}q^k\qbinom{2k}kq\eq0\ (\mo\ [3^a]_q^2)\quad\t{for any}\ m\in\Z^+.\tag1.3$$
\endproclaim

Concerning this conjecture, Guo and Zeng [GZ] were able to show (1.3) with the modulus $[3^a]_q^2$
replaced by $[3^a]_q$.

In this paper we confirm Conjecture 1.2 and give a $q$-analogue of (1.2).

 \proclaim{Theorem 1.3} Let $a\in\Z^+$. Then
 $(1.3)$ holds. Furthermore, we have the following $q$-analogue of $(1.2)$:
$$\f1{[3^a]_q^2}\sum_{k=0}^{3^a-1}q^k\Mm{2k}{k}q\eq 2R(a,q)\ \l(\mo\ \Phi_{3^a}(q)\r)\tag1.4$$
where
$$R(a,q):=\sum_{\Sb k=1\\ 3\mid k-1\endSb}^{3^a-1}q^{\frac{(k+2)(k-1)}{6}}\frac{(-1)^{k}}{[k]_q^2}
\bigg(1+\l(\frac{k-1}{3}-\frac{3^{a-1}+1}{2}\r)(1-q^k)\bigg).\tag1.5$$
\endproclaim

\Remark\ 1.4. Let $a\in\Z^+$. Then $\lim_{q\to1}R(a,q)\eq-1\ (\mo\
3)$ since
$$\sum\Sb k=1\\3\mid k-1\endSb^{3^a-1}\f{(-1)^k}{k^2}=\sum_{j=0}^{3^{a-1}-1}\f{(-1)^{3j+1}}{(3j+1)^2}
\eq-\sum_{j=0}^{3^{a-1}-1}(-1)^j=-1\ (\mo\ 3).$$ Also, for
$k\in\Z^+$ with $k\eq1\ (\mo\ 3)$, $[k]_q$ is relatively prime to
$[3^a]_q$ since $k$ is relatively prime to $3^a$. Therefore (1.4)
implies both (1.2) and (1.3) in the case $m=1$.

We are going to prove an auxiliary result in the next section and then show Theorem 1.3 in Section 3.

\heading{2. An auxiliary theorem}\endheading

\proclaim{Theorem 2.1} Let $a,m\in\Z^+$ and let $\psi:\Z\to\Z$ be a
function such that for any $k\in\Z$ and $j=1,\ldots,a$ we have
$$\psi(k)\equiv\psi(-k)\pmod{3^a}\ \ \t{and}\ \ \psi(k+3^j)\equiv\psi(k)\pmod{3^j}.$$ Then
$$
\sum_{k=1}^{3^am-1}q^{\psi(k)}\jacob{k}{3}\qbinom{2\cdot
3^am}{k}{q}\equiv 0\pmod{[3^a]_q^2}.\tag2.1
$$
In particular,
$$
\sum_{k=1}^{3^am-1}\jacob{k}{3}\qbinom{2\cdot3^am}{k}{q}\equiv 0\pmod{[3^a]_q^2}.\tag2.2$$
\endproclaim
\Proof. Clearly $[x]_q\equiv [y]_q\pmod{\Phi_d(q)}$ provided that
$x\equiv y\pmod{d}$. By the $q$-Lucas congruence (cf. [Sa]),
$$
\qbinom{x_1d+y_1}{x_2d+y_2}{q}\equiv
\binom{x_1}{x_2}\qbinom{y_1}{y_2}{q}\pmod{\Phi_d(q)}
$$
for $x_1,x_2,y_1,y_2\in\N$ with $0\leq y_1,y_2\leq d-1$. Recall that
$$
[3^a]_q=\prod_{j=1}^{a}\Phi_{3^j}(q).
$$
Since these $\Phi_{3^j}(q)$ are relatively prime and $[2\cdot
3^a m]_q\equiv 0\pmod{[3^a]_q}$, we only need to show that
$$
\sum_{k=1}^{3^am-1}\jacob{k}{3}\frac{q^{\psi(k)}}{[k]_q}\qbinom{2\cdot
3^am-1}{k-1}{q}\equiv 0\pmod{\Phi_{3^j}(q)}
$$
for every $j=1,\ldots,a$.

For any $1\leq j\leq a$ and $1\leq k\leq 3^am-1$ with $3\nmid k$,
write $k=3^j s+t$ where $1\leq t\leq 3^j-1$. Then, by the $q$-Lucas
congruence,
$$
\qbinom{2\cdot 3^am-1}{k-1}{q}\equiv\binom{2\cdot
3^{a-j}m-1}{s}\qbinom{3^j-1}{t-1}{q}\pmod{\Phi_{3^j}(q)}.
$$
And we have
$$\align&\qbinom{3^j-1}{t-1}{q}=\prod_{j=1}^{t-1}\frac{[3^j-j]_q}{[j]_q}
\\=&\prod_{j=1}^{t-1}\frac{q^{-j}([3^j]_q-[j]_q)}{[j]_q}\equiv(-1)^{t-1}q^{-\binom{t}2}\pmod{\Phi_{3^j}(q)}.
\endalign$$ Hence
$$
\align
&\sum_{k=1}^{3^am-1}\jacob{k}{3}\frac{q^{\psi(k)}}{[k]_q}\qbinom{2\cdot 3^am-1}{k-1}{q}\\
=&\sum_{s=0}^{3^{a-j}m-1}\sum_{t=1}^{3^j-1}\jacob{3^j s+t}{3}\frac{q^{\psi(3^j s+t)}}{[3^j s+t]_q}\qbinom{2\cdot 3^am-1}{3^j s+t-1}{q}\\
\equiv&\sum_{s=0}^{3^{a-j}m-1}\binom{2\cdot 3^{a-j}m-1}{s}
\sum_{t=1}^{3^j-1}\jacob{t}{3}\frac{(-1)^{t-1}q^{\psi(t)-\binom{t}2}}{[t]_q}\pmod{\Phi_{3^j}(q)}.
\endalign
$$
Clearly,
$$
\align
&2\sum_{t=1}^{3^j-1}\jacob{t}{3}\frac{(-1)^{t-1}q^{\psi(t)-\binom{t}2}}{[t]_q}\\
=&\sum_{t=1}^{3^j-1}\bigg(\jacob{t}{3}\frac{(-1)^{t-1}q^{\psi(t)-\binom{t}2}}{[t]_q}+\jacob{3^j-t}{3}\frac{(-1)^{3^j-t-1}q^{\psi(3^j-t)-\binom{3^j-t}2}}{[3^j-t]_q}\bigg)\\
\equiv&\sum\Sb t=1\\3\nmid t\endSb^{3^j-1}\jacob{t}{3}\bigg(\frac{(-1)^{t-1}q^{\psi(t)-\binom{t}2}}{[t]_q}+\frac{(-1)^{t-1}q^{\psi(t)-\binom{-t}2}}{-q^{-t}[t]_q}\bigg)
=0\pmod{\Phi_{3^j}(q)}.
\endalign
$$
So (2.1) holds.

Note that (2.2) is just (2.1) with $\psi$ replaced by the zero function from $\Z\to\Z$.
So (2.2) is also valid. This concludes the proof. \qed

\heading{3. Proof of Theorem 1.3}\endheading

\proclaim{Lemma 3.1} Suppose that $k\eq l\ (\mo\ 3^a)$ where
$k,l\in\Z$ and $a\in\Z^+$. Then
$$2k^2-k\l(\f k3\r)\eq 2l^2-l\l(\f l3\r)\ (\mo\ 3^{a+1}).$$
\endproclaim
\Proof. Observe that
$$2(k+l)-\l(\f k3\r)\eq 4k-\l(\f k3\r)\eq0\ (\mo\ 3).$$
Thus
$$\align&2k^2-k\l(\f k3\r)-\l(2l^2-l\l(\f l3\r)\r)
\\=&(k-l)\l(2(k+l)-\l(\f k3\r)\r)\eq0\ (\mo\ 3^{a+1}).
\endalign$$
We are done. \qed

\proclaim{Lemma 3.2} Let $a\in\Z^+$ and let $\psi$ be a function
as in Theorem 2.1. Then
$$
\aligned
&\f1{2[3^a]_q^2}\sum_{k=1}^{3^a-1}q^{\psi(k)}\jacob{k}{3}\qbinom{2\cdot 3^a}{k}{q}\\
\equiv&\sum\Sb k=1\\3\mid k-1\endSb^{3^a-1}
 q^{\psi(k)-\binom{k}{2}} \f{(-1)^{k-1}}{[k]_q^2}
(1+\Psi_a(k)(1-q^k)) \pmod{\Phi_{3^a}(q)},
\endaligned\tag3.1
$$
where
$$\Psi_a(k):=\frac{\psi(3^a-k)-\psi(k)}{3^a}+\frac{3^{a}-1}{2}-k.\tag3.2$$
\endproclaim
\Proof. We have
$$
\align
&\sum_{k=1}^{3^a-1}\jacob{k}{3}\frac{q^{\psi(k)}}{[k]_q}\qbinom{2\cdot 3^a-1}{k-1}{q}\\
=&\sum_{k=1}^{3^a-1}\jacob{k}{3}\frac{q^{\psi(k)}}{[k]_q}\prod_{j=1}^{k-1}\frac{q^{-j}([2\cdot 3^a]_q-[j]_q)}{[j]_q}\\
\equiv&\sum_{k=1}^{3^a-1}\jacob{k}{3}\frac{(-1)^{k-1}q^{\psi(k)-\binom{k}2}}{[k]_q}
\bigg(1-2\sum_{j=1}^{k-1}\frac{[3^a]_q}{[j]_q}\bigg)\pmod{\Phi_{3^a}(q)^2},
\endalign
$$
since
$$[2\cdot 3^a]_q=[3^a]_q(1+q^{3^a})=[3^a]_q(2+q^{3^a}-1)\equiv
2[3^a]_q\pmod{[3^a]_q^2}.
$$
Note that for $s=0,1,2,\ldots$ we have
$$\align
q^{3^as}=&1+(q^{3^a}-1)\sum_{j=0}^{s-1}q^{3^aj}=1+(q^{3^a}-1)\bigg(s+\sum_{j=0}^{s-1}(q^{3^aj}-1)\bigg)
\\\eq&1+s(q^{3^a}-1)\pmod{\Phi_{3^a}(q)^2}
\endalign$$ and
$$
q^{-3^as}\equiv\frac1{1+s(q^{3^a}-1)}=\frac{1-s(q^{3^a}-1)}{1-s^2(q^{3^a}-1)^2}\equiv1-s(q^{3^a}-1)\pmod{\Phi_{3^a}(q)^2}.
$$
Also, for each $1\leq k\leq 3^a-1$, we have
$$
\align
&\frac{q^{\psi(3^a-k)-\binom{3^a-k}2}}{[3^a-k]_q}\bigg(1-2\sum_{j=1}^{3^a-k-1}\frac{[3^a]_q}{[j]_q}\bigg)
\\=&\frac{q^{\psi(3^a-k)-\binom{3^a}{2}+3^a k-\binom{k+1}{2}}([3^a]_q+[k]_q)}{q^{-k}([3^a]_q^2-[k]_q^2)}
\bigg(1-2\sum_{j=k+1}^{3^a-1}\frac{[3^a]_q}{[3^a-j]_q}\bigg)\\
\equiv&\frac{q^{\psi(k)-\binom{k}{2}}(1+(\frac{\psi(3^a-k)-\psi(k)}{3^a}+k-\frac{3^a-1}{2})(q^{3^a}-1))([3^a]_q+[k]_q)}{-[k]_q^2}
\\&\times\bigg(1+2\sum_{j=k+1}^{3^a-1}\frac{q^j[3^a]_q}{[j]_q}\bigg)\\
\equiv&-\frac{q^{\psi(k)-\binom{k}{2}}}{[k]_q}\bigg(1+\l(\frac{\psi(3^a-k)-\psi(k)}{3^a}+k-\frac{3^a-1}{2}\r)(q^{3^a}-1)\bigg)
\\&-2\frac{q^{\psi(k)-\binom{k}{2}}}{[k]_q}\sum_{j=k+1}^{3^a-1}\frac{q^j[3^a]_q}{[j]_q}
-\frac{q^{\psi(k)-\binom{k}{2}}[3^a]_q}{[k]_q^2}
\pmod{\Phi_{3^a}(q)^2}.
\endalign
$$
Clearly,
$$
\sum_{j=k+1}^{3^a-1}\frac{q^j}{[j]_q}=\sum_{j=k+1}^{3^a-1}\frac{1+q^j-1}{[j]_q}=-(3^a-1-k)(1-q)+\sum_{j=k+1}^{3^a-1}\frac{1}{[j]_q},
$$
and
$$\align&
\sum_{j=1}^{3^a-1}\frac{1}{[j]_q}=\frac{1}{2}\sum_{j=1}^{3^a-1}\bigg(\frac{1}{[j]_q}+\frac{1}{[3^a-j]_q}\bigg)
\\\eq&\frac{1}{2}\sum_{j=1}^{3^a-1}\bigg(\frac{1}{[j]_q}-\frac{q^j}{[j]_q}\bigg)=\frac{3^a-1}{2}(1-q)\pmod{\Phi_{3^a}(q)}.
\endalign$$
Thus we get
$$
\align
&\frac{q^{\psi(k)-\binom{k}2}}{[k]_q}\bigg(1-\sum_{j=1}^{k-1}\frac{[2\cdot 3^a]_q}{[j]_q}\bigg)+\frac{q^{\psi(3^a-k)-\binom{3^a-k}2}}{[3^a-k]_q}\bigg(1-\sum_{j=1}^{3^a-k-1}\frac{[2\cdot 3^a]_q}{[j]_q}\bigg)\\
\equiv&-\frac{q^{\psi(k)-\binom{k}{2}}}{[k]_q}\l(\l(\frac{\psi(3^a-k)-\psi(k)}{3^a}+\frac{3}{2}(3^a-1)-k\r)(q^{3^a}-1)\r)
\\&-\frac{q^{\psi(k)-\binom{k}{2}}}{[k]_q}\l(2\sum_{j=1}^{k-1}\frac{[3^a]_q}{[j]_q}+2\sum_{j=k+1}^{3^a-1}\frac{[3^a]_q}{[j]_q}+\frac{q^{\psi(k)-\binom{k}{2}}[3^a]_q}{[k]_q^2}\r)\\
\equiv&\frac{q^{\psi(k)-\binom{k}{2}}[3^a]_q}{[k]_q^2}+\frac{q^{\psi(k)-\binom{k}{2}}}{[k]_q}\Psi_a(k)(1-q^{3^a})\pmod{\Phi_{3^a}(q)^2}.
\endalign
$$
It follows that
$$
\align
&\sum_{k=1}^{3^a-1}\jacob{k}{3}\frac{q^{\psi(k)}}{[k]_q}\qbinom{2\cdot 3^a-1}{k-1}{q}\\
\equiv&\sum\Sb k=1\\3\mid k-1\endSb^{3^a-1}(-1)^{k-1}\frac{q^{\psi(k)-\binom{k}{2}}}{[k]_q}\bigg(\frac{[3^a]_q}{[k]_q}+\Psi_a(k)(1-q^{3^a})\bigg)\pmod{\Phi_{3^a}(q)^2}.
\endalign
$$
Noting that $[3^a]_q$ divides both sides of the above congruence by
Theorem 2.1 and $[2\cdot 3^a]\equiv 2[3^a]\pmod{[3^a]_q^2}$, we are
done. \qed

\medskip
\noindent{\it Proof of Theorem 1.3}.
Let $m\in\Z^+$. By [T, (4.3)] in the case $d=0$, we have
$$
\sum_{k=0}^{3^am-1}q^k\qbinom{2k}{k}q=-\sum_{k=1}^{3^am-1}q^{\psi_m(k)}\jacob{k}{3}\qbinom{2\cdot
3^am}{k}{q},\tag3.3
$$
where
$$
\psi_m(k)=\frac{2(3^am-k)^2-(3^am-k)\jacob{3^am-k}{3}-1}3.
$$
According to Lemma 3.1, the function $\psi=\psi_m$ has the property
described in Theorem 2.1. Combining (2.1) with (3.3) we get (1.3).

Now it remains to prove (1.4). By (3.3) and Lemma 3.2, we finally
obtain
$$\align
&\sum_{k=0}^{3^a-1}q^k\qbinom{2k}{k}q=-\sum_{k=1}^{3^a-1}q^{\psi_1(k)}\jacob{k}{3}\qbinom{2\cdot 3^a}{k}{q}\\
\equiv&2[3^a]_q^2\sum_{\Sb k=1\\ 3\mid k-1\endSb}^{3^a-1}q^{\frac{(k+2)(k-1)}{6}}\frac{(-1)^{k}}{[k]_q^2}\bigg(1+\l(\frac{k-1}{3}-\frac{3^{a-1}+1}{2}\r)(1-q^k)\bigg)\\
&\pmod{\Phi_{3^a}(q)[3^a]_q^2}.
\endalign
$$
This concludes our proof. \qed
\medskip

\Ack. The authors are indebted to Prof. J. Shallit for informing the reference [SSZ].


\bigskip

 \widestnumber\key{S09b}

 \Refs

\ref\key GZ\by V. J. W. Guo and J. Zeng\paper Some congruences involving central $q$-binomial coefficients
\jour preprint, arXiv:0910.3563. {\tt http://arxiv.org/abs/0910.3563}\endref

\ref\key PS\by H. Pan and Z. W. Sun\paper A combinatorial identity
with application to Catalan numbers \jour Discrete Math.\vol
306\yr 2006\pages 1921--1940\endref

\ref\key Sa \by B. E. Sagan\paper Congruence properties of $q$-analogs
\jour Adv. in Math.\vol 95\yr 1992\pages 127--143
\endref

\ref\key SSZ\by N. Strauss, J. Shallit and D. Zagier
\paper Some strange $3$-adic identities\jour Amer. Math. MOnthly
\vol 99\yr 1992\pages 66--69\endref

\ref\key St\by R. P. Stanley\book Enumerative Combinatorics \publ
Vol. 2, Cambridge Univ. Press, Cambridge, 1999\endref

\ref\key S06\by Z. W. Sun\paper Binomial coefficients and quadratic fields
\jour Proc. Amer. Math. Soc.\vol 134\yr 2006\pages 2213--2222\endref

\ref\key S09a\by Z. W. Sun\paper Binomial coefficients, Catalan numbers and Lucas quotients
\jour preprint, arXiv:0909.5648. {\tt http://arxiv.org/abs/0909.5648}\endref

\ref\key S09b\by Z. W. Sun\paper Various congruences involving binomial coefficients and higher-order Catalan numbers
\jour preprint, arXiv:0909.3808. {\tt http://arxiv.org/abs/0909.3808}\endref

\ref\key S09c\by Z. W. Sun\paper $p$-adic valuations of some sums of multinomial coefficients
\jour preprint, arXiv:0910.3892. {\tt http://arxiv.org/abs/0910.3892}\endref

\ref\key ST1\by Z. W. Sun and R. Tauraso\paper On some new congruences for binomial coefficients
\jour Acta Arith.\pages to appear. {\tt http://arxiv.org/abs/0709.1665}\endref

\ref\key ST2\by Z. W. Sun and R. Tauraso\paper New congruences for central binomial coefficients
\jour Adv. in Math., to appear. {\tt http://arxiv.org/abs/0805.0563}\endref

\ref\key T\by  R. Tauraso\paper $q$-analogs of some congruences
involving Catalan numbers, prepint, arXiv:0905.3816. {\tt
http://arxiv.org/abs/0905.3816}\endref

\ref\key ZPS\by L. Zhao, H. Pan and Z. W. Sun\paper Some congruences for the second-order Catalan numbers
\jour Proc. Amer. Math. Soc.\vol 138\yr 2010\pages 37--46\endref

\endRefs
\enddocument